# A probabilistic approach of the Poincaré-Bendixon's problem in $\mathbb{R}^d$

guy Cirier[1]
*LSTA, University Pierre et Marie Curie Sorbonne, France*

**Abstract**: We present how a probabilistic model can describe the asymptotic behaviour of the iterations, especially for ODE with an approach of the Poincaré-Bendixon's problem in $\mathbb{R}^d$.

**Résumé:** On présente un modèle probabiliste pour décrire le comportement asymptotique d'une itération, en particulier pour les EDO et pour aborder le problème de Poincaré-Bendixon's dans $\mathbb{R}^d$.

**Introduction**

Let $f$ a function which applies a set $C \subset \mathbb{R}^d$. We iterate $f$ indefinitely. Sometimes, the process converges to some fixed point or to some cycle. But, in many cases, it is quite impossible to know, after a long time, the position of the iteration $f^{(n)}$.

One of the best deterministic method to study the asymptotic behavior of the iteration is to linearize *f.* It consists to find an invertible function $\varphi$ such as: $\varphi \circ f = \lambda\varphi$ where $\lambda$ is linear. When this linearization is possible, $\varphi \circ f^{(n)} = \lambda^n\varphi$ and we obtain asymptotic cases, which are the generalization of the unidimensional cases $|\lambda| >, =, < 1$. But, we meet some well-known difficulties: we have basin of attraction of each fixed point of $f$ with fuzzy frontier. If we have many fixed points, what is the good choice among them? If $\lambda^p = 1$, we have other difficulties called resonance in $\mathbb{R}^d$. It is often a good approach near each fixed point.

However, it is very important to know what happens when we iterate *f*, especially if the *set* $C \subset \mathbb{R}^d$ *is bounded and when the set* $C \subset \mathbb{R}^d$ *is applied in itself*. In this case, a probabilistic approach with invariant measures gives other information. This is the object of this paper.

**1 - The Perron-Frobenius's measure**

Let $P$ a measure on a bounded set $C \subset \mathbb{R}^d$ and $P_f = P \circ f^{-1}$ the transform of $P$ by the function $f$. We define the measure $P$ of Perron-Frobenius as:

*$P$ is invariant under $f$ if, for all borelian set $B$, $P$ verifies the Perron-Frobenius's equation (PF):*

$$P_f(B) = P \circ f^{-1}(B) = P(B)$$

This measure $P$ is invariant when we iterate the measurable function $f$ in $\mathbb{R}^d$. Under very general conditions, the solution of this equation is unique.

This measure presents the same difficulties as we have seen with linearization methods for $f$: we will see that it depends of the fixed points and meets the resonance's problems, but it gives us many information about the areas where the iteration belongs more frequently. This information is asymptotic when $n \to \infty$ but doesn't give us any result about the transient steps of null measure.

This invariant measure is generally difficult to study. For instance, in the very simple case where *f* is invertible, its density $p$ verifies the functional equation:

$$p = p_f = |f^{-1}| p \circ f^{-1}$$

This is very complicated to solve. if $f$ is not invertible, it is more difficult.

**The Fourier-Laplace's transform approach**

Here, we seek an analytic approach of the invariant measure with the Fourier-Laplace's transform. We use the known property of invariant measure under $f$ . For all positive $P$ - measurable function $g$, we have the well-known formula:

$$\int g \circ f(x)\, dP(x) = \int g(x)\, dP(x)$$

With $g(x) = e^{y(x)}$, we write the Fourier-Laplace's transform $\emptyset(y) = \mathcal{L}(e^{yX}) = E(e^{yX})$ with the series $\emptyset(y) = \Sigma_n b_n y^n$ and $\emptyset_f(y) = \mathcal{L}\left(e^{yf(X)}\right)$. If the measure is invariant:

$$\emptyset(y) = \emptyset_f(y)$$

**Hypothesis**

*All along the paper, we suppose that the set $C \subset \mathbb{R}^d$ is bounded and $f$ applies $C$ in $C$ and is at least indefinitely derivable.*

The series $\emptyset(y) = \Sigma_n b_n y^n$ is convergent because it is bounded by an exponential series with the diameter of $C$. We translate the distribution with a small fixed vector $a \epsilon \mathbb{R}^d$, $X \to X + a$.

So: $\quad \emptyset(y, a) = E\left(e^{y(X+a)}\right)$

And $\quad \emptyset(y, a) \to \emptyset_f(y, a) = E\left(e^{yf(X+a)}\right)$

If the measure is invariant: $\emptyset(y, a) = \emptyset_f(y, a)$

We study the problem near a fixed point 0 of $f$: $f(0) = 0$.

***Proposition:***

*The resolving equation $R_a$ of PF is: $\theta_f(y, a) = \emptyset(y, a) - \emptyset_f(y, a) = 0$ .*

*If $f$ is an indefinitely derivable function, then, for $\forall a \in C$ and $\forall y$:*

$$\theta_f(y, a) = \Sigma_n b_n \partial^n (e^{ya} - e^{yf(a)})/\partial a^n = 0$$

*Let $e^n(y, a) = \partial^n (e^{ya} - e^{yf(a)})/\partial a^n$ be the resolving gap*

$$\theta_f(y, a) = \Sigma_n b_n e^n(y, a) = 0$$

*At $a = 0$:* $\quad \theta_f(y) = \Sigma_n b_n e^n(y, 0) \equiv 0$

*And, as $\theta_f(y, a) \equiv 0$ is an identity, $\partial^p \partial^q \theta_f(y, a)/\partial y^p \partial a^q \equiv 0$ for p all and q.*

■ If the random variable $X \epsilon C \subset \mathbb{R}^d$ has a measure $P$ with density $p(x)$, the translated random variable $X + a$ has the same density. Using the convergent series: $\emptyset(y) = \Sigma_n b_n y^n$, we have for all small translation $a \in C \subset \mathbb{R}^d$ of the random vector $X$, the translated density:

$$p(x - a) = \mathcal{L}^{-1}\left(e^{ta}\, \emptyset(t)\right) = \mathcal{L}^{-1}(\, \Sigma_n b_n t^n e^{ta})$$

We have the distribution $\quad p(x - a) = \Sigma_n b_n \partial^n \delta(x - a)/\partial a^n$

As $\quad E\left(e^{y(X+a)}\right) = \Sigma_n b_n y^n e^{ya} = \Sigma_n b_n \partial^n e^{ya}/\partial a^n$

And $\quad E\left(e^{yf(X+a)}\right) = \int e^{yf(x+a)}\, dP(x) = \int e^{yf(x)}\, p(x - a)\, dx$

$$E\left(e^{yf(X+a)}\right) = \Sigma_n b_n \partial^n \left(\int e^{yf(x)}\, \delta(x - a) dx\right)/\partial a^n$$

$$E\left(e^{yf(X+a)}\right) = \Sigma_n b_n \partial^n e^{yf(a)}/\partial a^n$$

By difference, we get: $\quad \theta_f(y, a) = E\left(e^{y(X+a)}\right) - E\left(e^{yf(X+a)}\right)$ ■

**Remarks**

- We observe that $\partial^n (e^{yf(a)})/\partial a^n = H_n(y, a)\, e^{yf(a)}$ where $H_n(y, a)$ is a Bell-polynomial in $y$ with degree $n$. $e^n(y, 0) = y^n - H_n(y)$ is a polynomial with degree $n$. We can note $e^n(y, a) = \partial^n (e^{ya} - e^{yf(a)})/\partial a^n = y^n e^{ya} - H_n(y, a)\, e^{yf(a)}$ .

- We obtain $\emptyset_f(y)$ by putting $H_n(y)$ instead of $y^n$ in the series of $\emptyset(y) = \Sigma_n b_n y^n$.
- As $f(0) = 0$: $\emptyset_f(0) = \emptyset(0) = 1$, $b_0 = 1$. But the other $b_n$ are unknown.

**Consequence**
*The general solution of the linear equation $\theta_f(y) = 0$ has the form $b\varphi(y)$ where $b$ is an arbitrary constant real. So, we can write $\emptyset(y) = 1 + b\varphi(y)$ with an arbitrary constant $b$. It means that $\emptyset(y) = 1$, for all $\varphi(y) = 0$ . We have a lattice distribution for $\varphi(y) = 0$.*

The solution of the Perron-Frobenius's equation is particular case of the equation $\theta_f(y) = 0$. First, we study the effect of an iteration on $e^n(y, 0)$= 0 and on $\theta_f(y, 0)$.

***Proposition***
*Iteration $a_\ell \longmapsto f_\ell(a)$ acts as a derivation on $\theta_f(y, 0)$ and on $e^n(y,0)$*= 0 in the sense that :
$\theta_f(y, 0) \longmapsto \partial\theta_f(y, 0)/\partial a_\ell|_{a=0}$ and $e^n(y, 0) \longmapsto e^{n+1_\ell}(y, 0)$
$f^{(n)} \rightarrow f^{(n+1)}$induces $e^n(y, 0) \rightarrow e^{n+1}(y, 0)$ or $H_n(y, 0) \rightarrow H_{n+1}(y, 0)$
■The demonstration is easily based on the mean's formula for each coordinate for $a \rightarrow 0$. So, all the coordinates of *n* are equal.
For example, we study the impact of $a_\ell \longmapsto f_\ell(a)$ on $D = \theta_f\big(y, \overline{a_\ell}, f_\ell(a)\big)$.
As $\theta_f(y, \overline{a_\ell}, a_\ell) = 0$ :
$D = \theta_f\big(y, \overline{a_\ell}, f_\ell(a)\big)$-$\theta_f(y, \overline{a_\ell}, a_\ell)$
So: $D = (f_\ell(a) - a_\ell)\big(\partial\theta_f\big(y, \overline{a_\ell}, a_\ell + r((f_\ell(a) - a_\ell))\big)/\partial a_\ell\big)$
When $a \rightarrow 0$ $f_\ell(a) - a_\ell \sim a_\ell(\lambda_\ell - 1)$.
$D \sim a_\ell(\lambda_\ell - 1)\big(\partial\theta_f(y, 0)/\partial a_\ell\big)$
For similar raisons, if $e^n(y, a)$ = 0, then: $e^n(y, f_\ell(a)) \sim a_\ell(\lambda_\ell - 1)\partial e^n(y, 0))/\partial a_\ell$.
More, if we iterate $f$ , that means $a \longmapsto f(a)$ in $e^n(y, a)$, we obtain:
$$D \sim \Sigma_{\ell=1}^{\ell=d} a_\ell(\lambda_\ell - 1)\big(\partial\theta_f(y, 0)/\partial a_\ell\big)$$
If this quantity is null for all $a$, we must have *d* equations independent null.
We have the same result for $e^n(y, 0)$. ■

**2 - Solution de $R_0$**
We choose a sufficiently large index $n \in N^d$ , with $n = n_1$=...=$n_\ell$...= $n_d$

***Lemma***
*For a fixed $b_n \neq 0$, under non-resonance conditions, if a solution of $\theta *_n (y) = \Sigma_{m\leq n} b *_m e^m(y) = 0$ exists, then zeros of $e^n(y)$ are zeros of $\theta *_n (y)$ .*
■ The solution of this equation is obtained as the following:
We choose a sufficiently large index $n \in N^d$ such as: $\theta_{nf}(y) = \Sigma_{m\leq n} b_m e^m(y)$ verifies uniformly $|\theta_{nf}(y) - \theta_f(y)| < \epsilon$.
As $\theta_f(y) = 0$, we search an approximation $\theta *_n (y) = 0$, and estimators $b *_m$ such as we have:
$$\theta *_n (y) = \Sigma_{m\leq n} b *_m e^m(y) = 0$$
- For $y = 0$: $b_0 = 1$. As $\theta *_n (y)$ is a polynomial and, if we search a solution $\theta *_n (y) = 0$, ether all the coefficients of $\theta *_n (y)$ are null or the solution is valid only for the $y$ verifying $\theta *_n (y) = 0$. But, as the term of highest degree of $\theta *_n(y)$ is: $(1-\lambda^n)\, b_n y^n$ , we must have, under non-resonance conditions and for all small $b_n \neq 0$, $e^n(y) = 0$. (Because all the other gaps $e^m(y)$ have a lower degree for all $m < n$).
Then, zeros of $e^n(y)$ are zeros of $\theta *_n(y)$ . ■

***Theorem***

*Under the non-resonance condition, we can find a unique convergent solution of $\theta *_n(y) = \Sigma_{m\leq n} b *_m e^m(y) = 0$, up to an arbitrary constant b:*

$$\emptyset_n(y) = 1 - be^n(y)$$

*We obtain a lattice distribution defined by the zeros of $e^m(y)$*

*If $y^n\lambda^n \gg y^n$ :*

$$\emptyset_n(y) = 1 - bH_n(y).$$

*Then, the distribution of the real zeros of the polynomials $H_n(y)$ gets the distribution of the Perron-Frobenius's measure when n→ ∞.*

*We obtain a lattice distribution defined by the zeros of $H_n(y)$.*

*If $\lambda^n \gg 1$, real zeros of the polynomials $H_n(y)$ gets the distribution of the Perron-Frobenius's measure.*

■ We note the polynomials $\emptyset *_n(y) = 1 + \Sigma_{0<m\leq n} b *_m y^m$

and $\emptyset_f *_n(y) = 1 + \Sigma_{0<m\leq n} b *_m H_m(y)$

So: $\theta *_n(y) = \emptyset *_n(y) - \emptyset_f *_n(y)$

- We search a solution under the condition $e^n(y) = 0$

1- Now, for all $y$ verifying $e^n(y) = 0$ ,can we find a solution of $\theta *_n(y) = 0$ ?

And for all $m < n$, we note

$$A *_{n-1}(y) = \Sigma_{m<n} b *_m e^m(y) = \theta *_n(y) - b_n e^n(y).$$

If $\theta *_n(y) = 0$ :

$$A *_{n-1}(y) = \Sigma_{m<n} b *_m e^m(y) = -b_n e^n(y) = -(1-\lambda^n)\, b_n y^n - b_n \Sigma_{0<k<n} h_{nk}\, y^k$$

Where all the coefficients of $e^m(y) = y^m - H_m(y)$ are known because $H_m(y)$=$\Sigma_{0<k\leq m} h_{mk}\, y^k$.

So, we study in $A *_{n-1}(y)$ all the terms of $y^m$with degree $m< n$:

$$A *_{n-1}(y) = \Sigma_{m<n} b *_m (y^m - \Sigma_{0<k\leq m} h_{mk}\, y^k) = -b_n \Sigma_{0<k<n} h_{nk}\, y^k$$

For fixed arbitrarily $b *_n = b_n \neq 0$,

2- We obtain a finite triangular system of linear equations which can be solved step by step, and we can identify in a unique way all the unknown coefficients $b *_m$ in function of $b_n$ and the coefficients $h_{mk}$ of $H_m(y)$=$\Sigma_{0<k\leq m} h_{mk}\, y^k$ with $m\leq n \in N^d$ .

3 - This solution is unique for all $b *_n = b_n \neq 0$ arbitrarily fixed, near to the solution of $\theta_f(y) = 0$, as the $b *_m$ converge to the $b_m$. So, we can construct the polynomials $\emptyset *_n(y) - 1$ and $\emptyset_f *_n(y) - 1$ and we can write $\emptyset *_n(y) = 1 + b_n \varphi^n(y)$ where $b_n$ is arbitrary. That means $\emptyset *_n(y) = 1$ when $\varphi^n(y) = 0$; then, we can choose now $\varphi^n(y) = e^n(y)$.■

Different cases can happen according to $\lambda^n \gg 1$ or $\lambda^n \ll 1$ .

If all the coordinates of $|\lambda|$are less than 1, the process converges to the fixed point.

If some of them are less than 1, but others are greater than 1, we have an hyperbolic situation under no resonance conditions.

When $y^n\lambda^n \gg y^n$ , we can write for large $n$ :

$$\emptyset_n(y) \sim 1 - bH_n(y).\ ■$$

And now we have to study the zeros of $H_n(y)$.

**Remark** (see general reference)
We deduce, under general conditions, that, if $q(y)$ is the density of real zeros of $H_n(y)$ when $n\to\infty$, then the invariant density $p(x)$ of the Perron-Frobenius's measure is:

$$p(x) = (-x)\partial q(x)/\partial x$$

-For convenience of calculations, we suppose $f$ polynomial and $C \subset \mathbb{R}^d$ bounded.
The problem is reduced to find the asymptotic distribution of the zeros of $H_n(y)$.

**3 - Study of the zeros** of $H_n(y)$
We will see in the next paragraph, that all the zeros of $H_n(y)$ are distinct with the steepest descent's method applied to $H_n(y)$ and we get an estimation of the asymptotic distribution of the real zeros of $H_n(y)$.
- First, we use the steepest descent's method as Plancherel and Rotach use. We recall that the polynomial:

$$H_n(y) = e^{-yf(a)}\partial^n e^{yf(a)}/\partial a^n|_{a=0} = \partial^n e^{yf(a)}/\partial a^n|_{a=0}$$

can be represented by the Cauchy's integral:

$$H_{n-1}(y) = K\oint_\Gamma \frac{e^{yf(a)}}{a^n}da = K\oint_\Gamma e^{yf(a)-n\ln a}da$$

where $\Gamma$ is a closed polydisk around the fixed point 0 of $f$, $a \in \mathbb{C}^d$, $K$ can be taken as some finite non-null function, for all $n = (n_1,\ldots, n_\ell,\ldots, n_d)$. If the $n_\ell$ are not equal, we take $\mu = n_1+\ldots+n_\ell,\ldots+n_d$, and we fix: $z_\ell = n_\ell/\mu$. So, the integral can be written: $\oint_\Gamma e^{\mu\gamma(a)}\,da$ with the integrand called the Plancherel-Rotach's function, as:

$$\mu\gamma(a) = yf(a) - n\ln a = \Sigma_\ell\,(y_\ell a f_\ell(a) - n_\ell \ln a_\ell),$$

If $y_\ell = n_\ell s_\ell = \mu z_\ell s_\ell$ and $n_\ell \ln a_\ell = \mu z_\ell \ln a_\ell$, the Plancherel-Rotach's function is:

$$\mu\gamma(a) = \mu\Sigma_\ell\,(z_\ell s_\ell f_\ell(a) - z_\ell \ln a_\ell)$$

with:

$$\gamma(a) = \Sigma_\ell\,(z_\ell s_\ell f_\ell(a) - z_\ell \ln a_\ell) = \Sigma_\ell\, z_\ell \gamma_\ell(a)$$

which doesn't depend on $\mu$ because the $y_\ell$, thus $s_\ell$, can be taken arbitrarily. Here, if we choose a sufficiently large index $n \in N^d$, with $n = n_1=\ldots=n_\ell\ldots= n_d$, we have: $z_\ell = 1$
Then, the distribution of the $f$-invariant measure $P$ is given in general by the distribution of the zeros of $H_{n-1}(y)$ when $n\to\infty$.

- Second, we use the Plancherel-Rotach's method, which is the steepest descent's method. We search the critical point of $\gamma(a)$. Under the conditions of the general position, **t**he critical point $a_c$ maximizing $e^{n\gamma(a)}$ gives the solution. *The critical point* $a_c$ is defined by the equation:

$$\frac{\partial n\gamma(a)}{\partial a} = \frac{y\partial f(a)}{\partial a} - \frac{n}{a} = 0$$

(A sufficient condition to get this maximum is that the hessian matrix which is Hermitian of $\gamma(a)$ is definite negative at $a_c$). Let $s = y/n$ with $s_\ell = y_\ell/n$
The critical point must be isolated from the other critical points and at a finite distance.
We notice that the real part $\Re(\gamma)$ of $\gamma(a)$ cannot annul $H_{n-1}(y)$ and we have to conserve only the $\Re(\gamma)$ which make maximum the integrand. Only the imaginary part $\Im(\gamma)$ of $\gamma(a)$ can nullify $H_{n-1}(y)$:

$$\mu\Im(\gamma) = \Im(yf(a) - n\ln a) = k\pi$$

As each iteration $f_\ell$ acts as a derivation on $H_{n-1}(y)$, we see:

$$\Im(s_\ell a f_\ell(a) - \ln a_\ell) = \pi k_\ell/n$$

We obtain asymptotically $\Im(\gamma_\ell(a)) \to \kappa_\ell \pi$ when $n\to\infty$.
At the critical point $a$ where the $\kappa_\ell$ are identically independent uniform distribution on (0,1).

We can bond these distributions of the PF-equation to each fixed point $f(0) = 0$, then, we have local solutions. All these distributions can be masked by various situations. The principle of the maximum of the real part $\Re(\gamma)$ of $\gamma(a)$ provides a method to define the fuzzy frontiers of the different domains of attraction.
As everybody knows, the steepest descent's method is difficult to use, but it indicates all the varieties of behavior.
In the case of unidimensional function, the repartition of the zeros verifies:
$\Im(\gamma(a)) = \kappa\pi$ with $\partial\gamma/\partial a = 0$
So: $q(s)\, ds = Prob\,\{1\ zero \in (s, s + ds\} = d\kappa = \Im(d\gamma/ds)\, ds/\pi = \Im(f(a))\, ds/\pi$
Because: $d\gamma/ds = \partial\gamma/\partial s + \partial\gamma/\partial a.\partial a/\partial s = f(a)$

**4 - Examples:**

- Let the logistic map: $f(a) = \lambda a - a^2/2$; and $\gamma(a) = s\,(\lambda a - a^2/2) - \ln a$;
$\partial\gamma/\partial a = s\,(\lambda a - a^2) - 1 = 0$
we put $\lambda\sqrt{s} = 2cos\vartheta$ , we have: $2\cos\vartheta\, a\sqrt{s} - sa^2 - 1 = 0$ with roots: $a\sqrt{s} = e^{\pm i\vartheta}$
and: $\Im(f(a)) = \Im(\lambda a - a^2/2) = \sin 2\vartheta/s$
$q(s)\, ds = (1 - cos2\vartheta)d\vartheta/\pi$
So: $q(s) = ((\lambda/2\pi)\sqrt{1/s - \lambda^2/4}$
If we put $t = cos\vartheta = \frac{\lambda\sqrt{s}}{2}$ then $W(t) = (2/\pi)\,\sqrt{1 - t^2}\,dt$
We recover directly a well-known result: Let $H_n\,(y, a) = \partial^n (e^{y(\lambda a - a^2/2)})/\partial a^n$
where $e^{y(\lambda a - a^2/2)}$ is (with easy transformations) like the generatrix function $e^{(2ta - a^2)}$ of the Hermite polynomials $H_n\,(t)$. The law of the zeros of $H_n\,(x)$ is known as the semi-circular Wigner's law:

$$W(t)dt = (2/\pi)\sqrt{1 - t^2}dt$$

- Then, the density of the logistic corresponding to $q(s)$ is:

$$p(s) = -sdq/ds = -s(2/\pi)d(\frac{\lambda}{4\sqrt{s}}\sqrt{1 - \frac{s\lambda^2}{4}})/ds$$
$$p(s) = -s\lambda(-1/4s^2\pi)/2\sqrt{1/s - \lambda^2/4}$$
$$p(s) = \lambda/(2\pi\sqrt{4s - s^2\lambda^2})$$

We deduce that the density of the logistic map follows a Beta (1/2,1/2) low in a more general situation than in the Ulam-Von Neumann's case.
- If we have a map $a_1 = \lambda a + a^2/2$, the corresponding Hermite polynomials $H_n\,(x)$ are always positive except if $x = 0$.
- $m$ - Hermitian case: $f(a) = \lambda a - a^m/m$
The Plancherel-Rotach's function is: $\gamma(a) = s\,(\lambda a - a^m/m) - \ln a$
With the critical point $a$ defined by the trinomial equation $d\gamma(a)/da = s(\lambda a - a^m) - 1 = 0$ studied by H. Fell.

**Consequence**

We take now a function $f$ quadratic in $\mathbb{R}^d$ with $f(0) = 0$. We write the PR function $\gamma(a)$for every fixed point 0 of $f(a)$:

$$f(a) = \lambda a + Qa^2$$

the hessian of $sQ$ is symmetric. For all $s$ such as $sQ$ is non -degenerate, it exists an orthogonal transformation $T$: $a = Tu$, with: $T'sQT = D$ , the diagonal matrix of eigenvalues of $sQ$ and
$\ln a = \Sigma_{\ell=1}^{\ell=d} \ln a_\ell = \ln\Pi_{\ell=1}^{\ell=d} a_\ell = \ln\mathrm{Vol}(a) = \ln\mathrm{Vol(u)} = \ln\Pi_{\ell=1}^{\ell=d} u_\ell = \ln u$
Because the volume $\mathrm{Vol}(a) = \mathrm{Vol(u)}$ is invariant under an orthogonal transformation.

We note $D_\ell = K_\ell^2$ if $\ell = 1,..,p$ and $D_\ell = -K_\ell^2$ if $\ell = p+1,..,d$.
Then, the P.R. function $\gamma(a)$ becomes:

$$\gamma(u) = sf(Tu) - \ln Tu = s\lambda Tu + Du^2 - \ln u$$
$$= \Sigma_{\ell=1}^{\ell=d}\Lambda_\ell u_\ell + \Sigma_{\ell=1}^{\ell=p}K_\ell^2 {u_\ell}^2 - \Sigma_{\ell=p+1}^{\ell=d}K_\ell^2 {u_\ell}^2 - \Sigma_{\ell=1}^{\ell=d}\ln u_\ell$$
$$=\Sigma_{\ell=1}^{\ell=p}(\Lambda_\ell u_\ell + K_\ell^2 {u_\ell}^2\text{-}\ln u_\ell) + \Sigma_{\ell=p+1}^{\ell=d}(\Lambda_\ell u_\ell - K_\ell^2 {u_\ell}^2\text{-}\ln u_\ell)$$

Where $\Lambda u = s\lambda Tu.$

If we note: $\gamma_+(u_\ell) = \Lambda_\ell u_\ell + K_\ell^2 {u_\ell}^2\text{-}\ln u_\ell$

and $\gamma_-(u_\ell) = \Lambda_\ell u_\ell - K_\ell^2 {u_\ell}^2\text{-}\ln u_\ell$

So: $\gamma(u) = \Sigma_{\ell=1}^{\ell=p}\gamma_+(u_\ell) + \Sigma_{\ell=p+1}^{\ell=d}\gamma_-(u_\ell)$

And, applying the logistic calculus to each $\gamma_+(u_\ell)$ and $\gamma_-(u_\ell)$, we obtain *p* conditions $\Lambda_\ell u_\ell = 0$ and *d-p* random independent variables following a Beta (1/2,1/2) low. But, we may have other fixed points: $a(1-\lambda) = Qa^2$.

## 5 - A differential equation as an iteration

We consider ordinary differential equation:

$$da/dt = F(a)$$

where $a \in C \subset \mathbb{R}^d or\ \mathbb{C}^d,\ t \in \mathbb{R}^+,\ F(a)$ is a polynomial application of $a \in C$ in $C$. The domain $C$ is supposed bounded. The problem is to find a function $a(t)$ verifying this equation with an initial condition $a(t_0)= a_0$.

The solution is theoretically $a(t)$ for $t > t_0$ :

$$a(t)= a_0+\int_{t_0}^{t} F(a(u))du$$

### The differential iteration

We associate the differential iteration $f(a)$supposed belonging in the bounded domain $C$:

$$f(a) = a + \delta F(a)$$

where $\delta = t/n$ is the path. When we iterate n times, we have

$$a_n = f^{(n)}(a(t_0)) \text{ and } a_p = f^{(p)}(a(t_0))$$

The method gives the solution $a_n$ by iterating *n* times $f(a)$ from a starting point $a(t_0)$ with the path $\delta = t/n$ and this solution $a_n \to a(t)$ when $n \to \infty$:

For $n > p$ $\quad a_n = f^{(n)}(a_0) = a(t_0) + \delta\left(\Sigma_{p=0}^{p=n-1}F(a_p)\right)$

$$a_n = a(t_0) + \delta S_n(a_0)$$

Then, when $n \to \infty$ : $\quad a(t)= \lim_{n\to\infty} f^{(n)}(a_0)= a_0+\int_{t_0}^{t} F(a(u))du$

with $\quad \delta S_n(a_0) \to S(a(t), a_0) = \int_{t_0}^{t} F(a(u))du$

The fixed points of a differential iteration are the zeros $\alpha$ of $F$: $F(\alpha) = 0$

### The invariant measure of the differential iteration

Now, we submit a probabilistic version of the Poincaré-Bendixon's problem in $\mathbb{R}^d$.

### *Proposition*

*Under the previous hypothesis, all the non-null measures verify:*

$$E\left(\int_0^1 yS(X(t))\, e^{yX+vyS(X(t))}dv\right) = 0$$

*Then, we have asymptotic random cycles around each fixed points. For all these cycles, the times of return in each very small borelian set around a point of a cycle are constant in probability. Along each cycle, the conditional probability has a constant density.*

■ With $f(a)= a +\delta F(a)$ for every measurable function $F$. Then, for $a_n = f^{(n)}(a) = a(t_0) + \delta(\Sigma_{p=0}^{p=n-1} F(a_p))$ with $\delta = t/n$, we must have the resolving equation in the neighborhood each fixed point for one or $n$ iterations:

$$\emptyset(y) = E(e^{yX}) = E\left(e^{yf(X)}\right) = E\left(e^{y\,f^{(n)}(X)}\right)$$

That means, especially for $f^{(n)}$:

$$\theta_f(y) = E\left(e^{yX} - e^{y\,f^{(n)}(X)}\right) = E\left(e^{yX}(1 - e^{y\delta S_n(X)}\right) = 0$$

As: $$\delta S_n(a_0) \to S(a(t), a_0) = \int_{t_0}^{t} F\big(a(u)\big)du = S(X(t))$$

$$\theta_f(y) \to E\left(e^{yX}(1 - e^{yS(X(t))}\right) = E\left(\int_0^1 d(\,e^{yX+vyS(X(t))})/dvdv\right) = 0$$

But $$\theta_f(y) = E\left(e^{yX} - e^{y\,f^{(n)}(X)}\right) \to E\left(e^{yX} - e^{y(X(t))}\right) = 0$$

In consequence, if $E\left(\int_0^1 d(\,e^{yX+vyS(X(t))})/dvdv\right) = 0$, we have non-null measures verifying $\theta_f(y) = 0$. In other words, $a(t) = a(t_0)$ for the invariant measure and some $t_0$ and we have asymptotic random periodic cycles under this condition.

But, when we have many fixed points, the complete solution is more difficult because we meet some problems with domains of domination (See §3). ■

**Remark**

Theoretically, if we know the probability's measure, we can define the statistics (mean, standard deviation…).

## 6 -Examples

1- Suppose that $F$ has a hessian definite negative, then, when $\delta \to 0$, it is easy to verify that the critical point verifies , $ay = 1$, with an approximation of $H_n(y) = \partial^n yF(a)e^{ya+yv\delta f(a)}/\partial a^n|_{a=0}$ for $\delta v \to 0$ .

If $\delta$ =0, the critical point $a_c$ is real and we don't have probabilistic solution.

2- Suppose we have a linearity in $b$: $\boldsymbol{a} = (a, b)$. We write $d\,\boldsymbol{a}/dt = F(\boldsymbol{a})$ with $\boldsymbol{a} = (a, b)$ with $\boldsymbol{a} = (a, b)$ $a \in \mathbb{R}^d$ and $b \in \mathbb{R}$:

$$da/dt = f(a)b + g((a)$$
$$db/dt = h(a)b + k(a)$$

We write the Plancherel-Rotach's function with then $\boldsymbol{y} = (y, z), y \in \mathbb{R}^d$ and $z \in \mathbb{R}$

$$n\gamma(\boldsymbol{a}) = y(a + \delta(f(a)b + g\big((a)\big) + z(b + \delta(h(a)b + k(a))$$

Putting $z = \delta z'$ and $b' = \delta b$ for $\delta > 0$, we obtain when $\delta \to 0$:

$$n\gamma(\boldsymbol{a}') \to y(a + f(a)b') + z'b' \text{-} n \ln b' \text{-} n \ln a$$

(The change $b' = \delta b$ doesn't modify the equation $\theta_f(\boldsymbol{y}) = 0$ ). And the critical point is:

$$yf(a)b' + z'b' \text{-} n = 0$$
$$ya + yb'\partial f(a)/\partial a \text{-} n = 0$$

So: $$ya(yf(a) + z') + ny\partial f(a)/\partial a \text{-} n\,(yf(a) + z') = 0$$

The imaginary critical points give the distribution of the cycles. Under general conditions, this distribution doesn't depend on the functions $g((a), h(a), k(a)$ but only on $f(a)$.

## 7- Critical frequencies

Asymptotically, we have random cycles. Let $a(t)$ be a point on a such asymptotic cycle and a very small invariant borelian around this point. So, we have many large times to return in this borelian. In the differential iteration, we have many and large $\tau = (t + kT)/n$ which give the same $a(t)$ where is a random quasi-period.

***Proposition***
*When the number of iterations n→ ∞ and if the la hessian of yF is definite negative, the approximation with defines* $s = y/n$ *in function of the critical point a:*

$$s + \tau s\partial F(a)/\partial a - 1/a = 0$$

*where* $1/a = (1/a_\ell, \ell = 1,2,\ldots, d)$.

*If* $s_a$ *is a particular solution and if* $\vec{s}$ *is an eigenvector of* $-\partial F(a)/\partial a$ *for the eigenvalue* $1/t$*, the general solution is:*

$$s = s_a + \vec{s}$$

*The eigenvalue* $1/t$ *can be interpreted as a critical asymptotic frequency.*
■ Contrary to the previous §5, we don't write the critical point $a$ as a function of $s$, but $s$ as a function of $a$. For fixed $a$ on an asymptotic cycle, we recognize a linear affine equation of $s$ depending on the parameter $t$. We have to find a particular solution $s_a$ :

$$s_a + \tau s_a\partial F(a)/\partial a - 1/a = 0$$

Formally: $$s_a = (Id + \tau\partial F(a)/\partial a)^{-1}1/a$$

The equation of $s_a$ is now an elementary equation and has a unique solution for all $t \neq -1/\lambda_a$ where $\lambda_a$ is eigenvalue of $\partial F(a)/\partial a$ at the fixed critical point $a$,
Let the general solution be $s = s_a + \vec{s}$, where $s_a$ is a particular solution of the general equation:

As: $$s_a + \vec{s} + \tau(s_a + \vec{s})\partial F(a)/\partial a - 1/a = 0$$

then : $$\vec{s} = -\tau\vec{s}\partial F(a)/\partial a$$

$\vec{s}$ is eigenvector of $\partial F(a)/\partial a$ for the eigenvalue $\lambda_a$ and giving $\tau = -1/\lambda_a$ maximal positive; $\vec{s}$ is defined with a multiplicative constant arbitrary. The general solution is: $s = s_a + \vec{s}$ and shows a discontinuity at the eigenvalues $\lambda_a$.■

**Remark**: calculation of $s_a$
$s_a$ is obtained with $(Id + \tau\partial F(a)/\partial a)^{-1}$ for all $\tau \neq -1/\lambda_a$ which doesn't belong to the spectrum of $-\partial F(a)/\partial a$ with the series development of $\tau$.

## 8- Case where the hessian is degenerate: the equation of Lorenz

Generally, the hessian is not definite negative. The Lorenz's equation is an example particularly important because the differential iteration can be broken down into three independent iterations which have a remarkable feature: a partial linearity; an iteration with a negative hessian which induces a probabilistic solution and another with a positive hessian. It is an ideal example to clarify the previous results.
However, as there is an interpenetration of the distributions related to each fixed point, the connection between the various results remains delicate. The probabilistic presentation seems to be the least bad: it gives the probability of presence except at the places where the domination changes; in this case, we go from a basin of iteration to an another.

*- Presentation of the differential iteration at its fixed points.* These equations are written in bold notations:

$$d\,\boldsymbol{a}/dt = F(\boldsymbol{a}) \text{ where } \boldsymbol{a} = (a, b, c):$$
$$da/dt = \sigma(b\text{-}a)$$
$$db/dt = \rho a - b\text{-}ac$$
$$dc/dt = \text{-}\beta c + ab$$

*the differential iteration* $\boldsymbol{a}_1 = f(\boldsymbol{a})$ associated with a given path $\delta = t/n$ is:

$$a_1 = a + \delta\sigma(b\text{-}a)$$
$$b_1 = b + \delta(\rho a - b\text{-}ac)$$

$$c_1 = c + \delta(-\beta c + ab)$$

We recall the known results concerning the fixed points:

The fixed points are zeros of $F(\boldsymbol{a})=0$. If $\rho > 1$ and $\alpha = \sqrt{\beta(\rho-1)}$ , it exists three fixed points:

The point $O = (0,0,0)$, and two others symmetric to the axis of $c$:

$\alpha_+ = (\alpha, \alpha, \alpha^2/\beta)$ et $\alpha_- = (-\alpha, -\alpha, \alpha^2/\beta)$.

At $O$, the eigenvalue's equation $\lambda$ of the linear part is:

$(\beta + \lambda)\ [(\sigma + \lambda)\ (1 + \lambda) - \sigma\rho] = 0$,

But, at $\alpha_+$ or at $\alpha_-$:

$$\lambda(\beta + \lambda)(1 + \sigma + \lambda) - \alpha^2(2\sigma + \lambda) = 0,$$

Coefficients $\beta, \sigma, \rho$ are such as these three repellent fixed points; that means we have to study the distributions around each fixed point. We don't speak here about attractive cycles, resonances, and some particular values of the parameters, etc. It remains many things to clarify.

The iteration applies a compact set C in itself for $\delta > t > 0$ (the phenomenon occurring between a cold sphere at -50° and hot sphere, the earth, at +15° as the terrestrial atmosphere is modelled).

This iteration is quadratic, but has a linearity in $\boldsymbol{a}$.

*Analysis of the hessian*

Projecting $f(\boldsymbol{a})$ onto an axis $\boldsymbol{y} = (x, y, z)$, we write:

$\boldsymbol{y} f(\boldsymbol{a}) = L(\boldsymbol{a}) + \delta Q(\boldsymbol{a})$

where $L(\boldsymbol{a})$ is linear for $\boldsymbol{a}$: $L(\boldsymbol{a}) = x(a + \delta\sigma(b-a)) + y(b + \delta(\rho a - b)) + zc(1 - \delta\beta)$

$L(\boldsymbol{a}) = aL_1 + bL_2 + cL_3$

with : $L_1 = x(1-\delta\sigma) + \delta\rho y$

$L_2 = \delta\sigma x + y(1-\delta)$

$L_3 = z\,(1-\delta\beta)$

and $Q(\boldsymbol{a})$ is quadratic: $\delta Q(\boldsymbol{a}) = \delta(zb - yc)\ a$

The hessian is degenerated and not definite negative. We cannot apply the previous results. On the other hand, we can always use the lemma which requires: $\partial(\partial^n(e^{yf(a)})/\partial a^n)\partial\delta|_{a=0} = 0$.

**A-** Before to study this equation, we examine the quadratic application and its matrix $Q(\boldsymbol{a})$ :

$$Q = \begin{bmatrix} 0 & z & -y \\ z & 0 & 0 \\ -y & 0 & 0 \end{bmatrix}$$

If $\mu = \sqrt{y^2 + z^2}$ is the positive eigenvalue of the characteristic equation of $Q$: $\mu\,(\mu^2 - y^2 - z^2) = 0$

The matrix of the eigenvectors $T$ is orthogonal and constant for all $\boldsymbol{a}$.

$$T = \frac{1}{\mu\sqrt{2}} \begin{bmatrix} 0 & \mu & \mu \\ y\sqrt{2} & -z & z \\ z\sqrt{2} & y & -y \end{bmatrix}$$

Corresponding to the diagonal matrix of the eigenvectors.

$$\Lambda = \begin{bmatrix} 0 & 0 & 0 \\ 0 & -\mu & 0 \\ 0 & 0 & \mu \end{bmatrix}$$

*- Changing of basis*

We calculate directly with the Hermite's polynomials.

The application $\boldsymbol{u} = T\boldsymbol{a}$ with $\boldsymbol{u} = (u, v, w)$ is orthogonal and transforms:
- $\boldsymbol{yf}(\boldsymbol{a})$ in: $G(\boldsymbol{u}) = \boldsymbol{y}f(T'\boldsymbol{u})$ :
- $Q(\boldsymbol{a})$ in: $Q(\boldsymbol{u}) = \delta\mu(w^2 - v^2)$
- $\boldsymbol{L}(\boldsymbol{a})$ in: $\boldsymbol{LT'u}$

Where $T'$ means the transposed of the orthogonal matrix $T$, which is also its inverse: $T' = T^{-1}$.
Now, in the new basis $\boldsymbol{u}$, the function $\boldsymbol{y}f(\boldsymbol{a})$ is factorized: $\boldsymbol{y}f(T'\boldsymbol{u}) = G(\boldsymbol{u})$ into three independent functions:

$$G(\boldsymbol{u}) = g_1(u) + g_2(v) + g_3(w)$$

with:

$$g_1(u) = l_1 u$$
$$g_2(v) = l_2 v - \delta\mu\, v^2$$
$$g_3(w) = l_3 w + \delta\,\mu\, w^2$$

We get 3 independent iterations:
- the first is linear;
- the second is a random iteration;
- the third remains positive, except if $l_3 = 0$.

To calculate $l_1$, $l_2$ et $l_3$, we form: $L(\boldsymbol{a}) = a\,(x - \delta\sigma x + \delta\rho y) + b(\delta\sigma x + y(1-\delta)) + zc(1-\delta\beta)$
With: $L_1 = x\,(1-\delta\sigma) + \delta\rho y$ ; $L_2 = \delta\sigma x + y(1-\delta)$ ; $L_3 = z(1-\delta\beta)$

Then: $$lu = (l_1, l_2, l_3)\boldsymbol{u} = LT'\boldsymbol{u} = (L_1, L_2, L_3)\frac{1}{\mu\sqrt{2}}\begin{bmatrix} 0 & y\sqrt{2} & z\sqrt{2} \\ \mu & -z & y \\ \mu & z & -y \end{bmatrix}$$

$$l_1 = (\delta\sigma x + y\,(1-\delta) + z(1-\delta\beta))/\sqrt{2}$$
$$l_2 = (x - \delta\sigma x + \delta\rho y)y/\mu - (\delta\sigma x + y(1-\delta)) - z(1-\delta\beta))z/\mu\sqrt{2}$$
$$l_3 = (x - \delta\sigma x + \delta\rho y)z/\mu + (\delta\sigma x + y(1-\delta)) - z(1-\delta\beta))y/\mu\sqrt{2}$$

**B-** Let the resolving gap $e^n(\boldsymbol{y}) = \partial(\partial^n(e^{yf(\boldsymbol{a})})/\partial\boldsymbol{a}^n)\partial\delta|_{a=0} = 0$
For $\forall\, t \le \delta$. Putting $\boldsymbol{a} = T'\boldsymbol{u}$, we have:

$$e^n(\boldsymbol{u}) = T^n\partial(\partial^n(e^{yf(T'\boldsymbol{u})})/\partial\boldsymbol{u}^n)\partial\delta|_{\boldsymbol{u}=0} = 0$$

$$\partial^n(e^{yf(T'\boldsymbol{u})})/\partial\boldsymbol{u}^n = \partial^n(e^{g_1(u)})/\partial u^n \,.\, \partial^n(e^{g_2(v)})/\partial v^n .\ \partial^n(e^{g_3(w)})/\partial w^n$$

This gives: $$\partial^n(e^{g_1(u)})/\partial u^n = {l_1}^n e^{g_1(u)}\ ;$$
$$\partial^n(e^{g_2(v)})/\partial v^n = H_n(g_2(v))e^{g_2(v)}\ ;$$
$$\partial^n(e^{g_3(w)})/\partial w^n = H_n(g_3(w))e^{g_3(w)}$$

And: $$e^n(\boldsymbol{u}) = \partial\, {l_1}^n H_n(g_2(v))H_n(g_3(w))(e^{yf(\boldsymbol{T'u})})\partial\delta|_{u=0} = 0$$

***Proposition***
*The solution around the fixed point 0 consists of the intersection of the family of random surfaces defined by: $l_2/2\sqrt{\mu} \longmapsto$ low $\beta(1/2,1/2)$ with The surfaces $\sigma x - y - z\beta = 0$ et $(-\sigma x + \rho y)z + (\sigma x - y + z\beta)y/\sqrt{2} = 0$.*

■ With the same calculations of encodings and interchanging the derivations, we have:

$$\partial {l_1}^n/\partial\delta = 0;\ \partial\, H_n(g_2(v))/\partial\delta = 0;\ \partial\, H_n(g_3(w))/\partial\delta = 0$$

We study separately the three expressions:
- First: $\partial\, {l_1}^n/\partial\delta = n(\partial l_1/\partial\delta){l_1}^{n-1} = 0$

Either $\partial l_1/\partial\delta = \sigma x - y - z\beta = 0$, or: $l_1 \sim (y + z)/\sqrt{2} = 0$
- Second: the polynomial $H_n(g_3(w))$ when $w = 0$ is a Hermite's polynomial $H_n(x)$ where $x$ is $x = il_3/(\sqrt{2\delta\mu}$ . this polynomial $i^n H_n(il_3/(\sqrt{2\delta\mu})$ is always positive whatever $n$. In a general way:

$$\partial\, H_n(x)/\partial\delta = nH_{n-1}(x)\ \partial\, x/\partial\delta = 0.\ \text{So:}\ d(l_3/\sqrt{2\delta\mu})/d\delta = 0,$$

And $l_3 \sim (xz\sqrt{2} + (y\text{-}z)\, y)/\mu \sqrt{2}=0$

- Third: in the case of $H_n(g_2(w))$, in addition to the solution $l_2 = 0$, we have to find the possible invariant distribution of $H_n(l_2/\sqrt{2\delta\mu}) = 0$.

Let the integrand of $n\gamma(w) = g_2(w) - n\ln w$

When $\delta \to 0$, $l_2 \sim (xy\sqrt{2} + (y\text{-}z)\,z)/\sqrt{2}\mu$ with $\mu = \sqrt{y^2 + z^2}$.

By normalization of the coordinates $\boldsymbol{x} = (x, y, z) = \delta n\boldsymbol{s} = (\delta nr, \delta ns, \delta nt)$, we obtain:

$$l_2 \sim n\delta(rs\sqrt{2} + (s - t)t\,)/\,2(s^2 + t^2)^{\frac{1}{2}} = n\delta l_2(\boldsymbol{s})$$

$$\delta\mu = n\delta^2(s^2 + t^2)^{1/2} = n\delta^2\,\mu(\boldsymbol{s})$$

$$n\gamma(v) = n(\delta l_2(\boldsymbol{s})v - \mu(\boldsymbol{s})(\delta v)^2 - \ln\delta v + \ln\delta)$$

Putting $\delta v = \mathrm{v}$, we have: $n\gamma(\mathrm{v}) = n(l_2(\boldsymbol{s})\mathrm{v} - \mu(\boldsymbol{s})\mathrm{v}^2 - \ln\mathrm{v})$

We search the critical point: $d\gamma(\mathrm{v})/d\mathrm{v} = l_2(\boldsymbol{s}) - 2\mu(\boldsymbol{s})\mathrm{v} - 1/\mathrm{v} = 0$

The imaginary roots are: $\mathrm{v}(\boldsymbol{s}) = l_2(\boldsymbol{s})/4\mu(\boldsymbol{s}) \pm i\sqrt{1/2\mu(\boldsymbol{s}) - l_2(\boldsymbol{s})^2/16\mu(\boldsymbol{s})^2}$

Under the condition: $l_2(\boldsymbol{s})^2 < 8\mu(\boldsymbol{s})$:

$$l_3 \sim (rt\sqrt{2} + (s\text{-}t)\,s)/\mu(\boldsymbol{s})\sqrt{2}=0$$

Implies: $l_2(\boldsymbol{s}) = -(s - t)^2\,/\,\sqrt{2}(s^2 + t^2)^{1/2}$

The condition becomes: $(s - t)^4\,/\,(s^2 + t^2)^{3/2} < 16$

$l_1 = 0$ implies $s + t = 0$, then: $s < 8$

In any case, we observe that the conditions $l_3 = l_1 = 0$ allow us to express $r$ et $t$ depending on $s$ and we can write that the density of zeros of $s$ is now:

$$q(s)ds = \text{Prob } (1 \textit{ zero between } s, s\text{+d}s) = |Imf(\mathrm{v}(s))|\mathrm{d}s/\pi$$

$$q(s)ds = l_2(\boldsymbol{s})\sqrt{8\mu(\boldsymbol{s}) - l_2(\boldsymbol{s})^2}/8\pi\mu(\boldsymbol{s}) = \mathrm{d}\kappa$$

$\kappa$ follows a uniform low on (0,1) with: $s + t = 0$ (or $\sigma x - y\text{-}z\beta = 0$) and: $xy\sqrt{2} + (y\text{-}z)z + 0$

We also remark that the normalization doesn't affect the coefficients of the orthogonal matrix:

$$T(x, y, z) = T(\delta nr, n\delta s, n\delta\, t) = T(r, s, t)$$

We now verify similar results the two other fixed points $\boldsymbol{\alpha}_+$ et $\boldsymbol{\alpha}_-$. ■

**C** - Calculation for the two other fixed points

We search the distributions around the two other fixed points. To pass from the fixed point $O$ to the fixed point $\alpha_+$ or $\alpha_-$, it is sufficient to put in the differential iteration instead of $\boldsymbol{a} = (a,b,c)$

$\boldsymbol{a}'\text{+}\boldsymbol{\alpha}_+ = (a'\text{+}\alpha,\ b'\text{+}\alpha,\ c'\text{+}\alpha^2/\beta)$ and $\boldsymbol{a}'\text{+}\boldsymbol{\alpha}_- = (a'\text{-}\alpha,\ b'\text{-}\alpha,\ c'\text{+}\alpha^2/\beta)$ :

So, for $\boldsymbol{a}'\text{+}\boldsymbol{\alpha}_+$ : $\boldsymbol{a}_1 = f(\boldsymbol{a})$ where $\boldsymbol{a}_1 = (a_1, b_1\, c_1)$ becomes

$\boldsymbol{a}_1 = \boldsymbol{a}'_1 + \alpha_+ = f(\boldsymbol{a}) = f(\boldsymbol{a}'\text{+}\alpha_+)$;

then : $\boldsymbol{a}'_1 = a'\text{+}\delta F(\boldsymbol{a}'\text{+}\alpha_+)$

And $\boldsymbol{a}_1 = f(\boldsymbol{a})$: $a_1 = a + \delta\sigma(b\text{-}\,a)$

$$b_1 = b\text{+}\ \delta(\rho a -\ b\text{-}ac)$$

$$c_1 = c\text{+}\ \delta(\text{-}\beta c + a\ b)$$

Becomes for $\boldsymbol{a}\ \text{+}\alpha_+$ (we remove apostrophe of $\boldsymbol{a}'$ to make the notation less cluttered):

$$a'_1 = a + \delta\sigma(b\text{-}\,a) = a_1$$

$$b'_1 = b\text{+}\ \delta(\rho a -\ b\ \text{-}\,a\ c) + \delta(\text{-}\,\alpha c\ \text{-}a\alpha^2/\beta) = b_1 + \delta(\text{-}\,\alpha c\ \text{-}a\alpha^2/\beta)$$

$$c'_1 = c\text{+}\ \delta(\text{-}\beta c\text{+}ab) + \delta\alpha\ (a\text{+}b) = c_1\text{+}\ \delta\alpha\ (a\text{+}b)$$

The projection of $f\ (\boldsymbol{a})$ on an axis $\boldsymbol{y} = (x, y, z)$ can be written:

$$\boldsymbol{y}f(\boldsymbol{a}') = xa_1\text{+}yb_1 + \delta y(\text{-}\,\alpha c\ \text{-}a\alpha^2/\beta)\text{+}zc_1\text{+}\ z\delta\alpha(a\text{+}b)$$

$$\boldsymbol{y}f\ (\boldsymbol{a}') = \boldsymbol{y}f(\boldsymbol{a})\text{+}\ \delta(\mathrm{a}\ (z\alpha\text{-}y\alpha^2/\beta)\ \text{+}z\alpha b\text{-}y\alpha c)$$

and $Q(\boldsymbol{a})$ is invariant: $\boldsymbol{y}f\boldsymbol{a}') = L'(\boldsymbol{a}) + \delta Q(\boldsymbol{a})$

$L(\boldsymbol{a})$ is linear for $\boldsymbol{a}$: $L'(\boldsymbol{a}) = L(\boldsymbol{a})\text{+}\ \delta(\mathrm{a}\ (z\alpha\text{-}y\alpha^2/\beta)\ \text{+}z\alpha b\text{-}y\alpha c)$

$$L'(\boldsymbol{a}) = aL'_1 + bL'_2 + cL'_3$$

with: $L'_1 = L_1 + \delta(z\alpha - y\alpha^2/\beta)$; $L'_2 = L_2 + \delta z\alpha$; $L'_3 = L_3 - \delta y\alpha$

then $T$ and $\Lambda$ remain invariant.

We calculate $l'_1, l'_2$ et $l'_3$, with $L(\boldsymbol{a}) = a(x - \delta\sigma x + \delta\rho y) + b(\delta\sigma x + y(1-\delta)) + zc(1-\delta\beta)$ :

Where $L_1 = x(1-\delta\sigma) + \delta\rho y$ ; $L_2 = \delta\sigma x + y(1-\delta)$ ; $L_3 = z(1-\delta\beta)$

And: $\boldsymbol{l'}u = (l'_1, l'_2 \text{ et } l'_3)\boldsymbol{u}$

$= LT'\boldsymbol{u}$

$$= (L_1 + \delta(z\alpha - y\alpha^2/\beta), L_2 + \delta z\alpha,\ L_3 - \delta y\alpha)\frac{1}{\mu\sqrt{2}}\begin{bmatrix} 0 & y\sqrt{2} & z\sqrt{2} \\ \mu & -z & y \\ \mu & z & -y \end{bmatrix}$$

The results are modified:

$$l'_1 = l'_1 + \delta\alpha(z - y)/\sqrt{2}$$
$$l'_2 = l'_2 + \delta\alpha((z - y\alpha/\beta)y\sqrt{2} - z(z + y))/\mu\sqrt{2}$$
$$l'_3 = l_3 + \delta\alpha((z - y\alpha/\beta)z\sqrt{2} + y(z + y))/\mu\sqrt{2}$$

The following calculations remains the same.

When $\boldsymbol{a}$ becomes $\boldsymbol{a} + \alpha-$ $a''_1 = a + \delta\sigma(b - a) = a_1$

$$b''_1 = b + \delta(\rho a - b - ac) + \delta(\alpha c - a\alpha^2/\beta) = b_1 + \delta(\alpha c - a\alpha^2/\beta)$$
$$c''_1 = c + \delta(-\beta c + ab) - \delta\alpha\ (a+b) = c_1 - \delta\alpha\ (a+b)$$

It remains the problems of domination and frontiers between the various distributions attached at each fixed point.

**Remark**

We have to go back to the original coordinates. And the solution gives only probabilities of presence...

## 9 - Conclusion

With this probabilistic method, we get many new results, but also many new difficulties such as, for instance, when the hessian is degenerated. The results located around each fixed point are other difficulties, but the steepest descent gives us the basin of each fixed point. And we don't must forget that the steepest descent's method is not always easy to use.

We have applied these results to PDE equations and obtain other new results.